# A note on $B_k$ − *sequences*

Li An-Ping

Beijing 100085, P.R. China
apli0001@sina.com

Abstract

A sequence of non-negative integers is called a $B_k$ − *sequences* if all the sums of arbitrary $k$ elements are different. In this paper, we will present a new estimation for the upper bound of $B_k$ − *sequences*.

# 1. Introduction

A sequence $\mathcal{A}$ is called $B_k - sequence$ if all the sums of arbitrary $k$ elements of $\mathcal{A}$ are different. Suppose $n$ is a arbitrary positive integer, denoted by $\Phi_k(n)$ the maximum of the sizes of the $B_k - sequences$ contained in $[0, n]$.

For the lower bound of $\Phi_k(n)$, J. Singer [10] and S.C. Bose and S. Chowla [1] proved

$$\Phi_k(n) \geq n^{1/k} + o(n^{1/k}). \tag{1.1}$$

For the upper bound of $\Phi_k(n)$, Erdos, Turan [4] provided that

$$\Phi_2(n) \leq n^{1/2} + O(n^{1/4}). \tag{1.2}$$

B. Lindstrom [8], [9] shown that

$$\Phi_2(n) \leq n^{1/2} + n^{1/4} + 1, \tag{1.3}$$

$$\Phi_4(n) \leq (8n)^{1/4} + O(n^{1/8}). \tag{1.4}$$

In [6], we presented

$$\Phi_3(n) \leq \left(\left(1 - \frac{1}{6\log_2^2 n}\right) 4n\right)^{1/3} + 7, \tag{1.5}$$

lately we improved that [7]

$$\Phi_3(n) \leq (3.962n)^{1/3} + 1.5. \tag{1.6}$$

J. Cilleruelo [3] presented that

$$\Phi_3(n) \leq \left(\frac{4n}{1 + 16/(\pi+2)^4}\right)^{1/3} + o(n^{1/3}) \tag{1.7}$$

$$\Phi_4(n) \leq \left(\frac{8n}{1 + 16/(\pi+2)^4}\right)^{1/4} + o(n^{1/4}), \tag{1.8}$$

For general $k$, X.D. Jia [5] and S. Chen [2] proved that

$$\Phi_{2m}(n) \leq \left(m \cdot (m!)^2\right)^{1/2m} \cdot n^{1/2m} + O(n^{1/4m}), \tag{1.9}$$

$$\Phi_{2m-1}(n) \leq \left((m!)^2\right)^{1/(2m-1)} \cdot n^{1/(2m-1)} + O(n^{1/(4m-2)}), \tag{1.10}$$

J. Cilleruelo [3] gave an improvement

$$\Phi_{2m}(n) \leq \begin{cases} \left(\dfrac{m \cdot (m!)^2}{1+\cos^{2m}(\pi/m)}\right)^{1/2m} \cdot n^{1/2m} + o(n^{1/2m}), & \text{for } m < 38 \\ \left(\dfrac{5}{2}\left(\dfrac{15}{4} - \dfrac{5}{4m}\right)^{1/4} \sqrt{m} \cdot (m!)^2 \right)^{1/2m} \cdot n^{1/2m} + o(n^{1/2m}) & \text{for } m \geq 38 \end{cases} \quad (1.11)$$

$$\Phi_{2m-1}(n) \leq \begin{cases} \left(\dfrac{(m!)^2}{1+\cos^{2m}(\pi/m)}\right)^{1/(2m-1)} \cdot n^{1/(2m-1)} + o(n^{1/(2m-1)}), & \text{for } m < 38 \\ \left(\dfrac{5}{2}\left(\dfrac{15}{4} - \dfrac{5}{4m}\right)^{1/4} \cdot (m!)^2 \right)^{1/(2m-1)} \cdot n^{1/(2m-1)} + o(n^{1/(2m-1)}) & \text{for } m \geq 38 \end{cases} \quad (1.12)$$

In this paper, we will present a new estimation for the upper bound of $\Phi_k(n)$ with some slightly different means.

## 2. The Main Results

At first, we introduce some notations will be used in this paper.

**Definition 1.** Suppose that $m$ and $\delta$ are two positive integers, a positive integer array $(\alpha_1, \alpha_2, \cdots, \alpha_i, \cdots)$ is called $\delta$-*deviation* set of size $m$, if it satisfying following equation

$$\begin{cases} \sum_i \alpha_i = m, \\ \left|\sum_i \alpha_{2i-1} - \sum_i \alpha_{2i}\right| \leq 1, \\ \left|\sum_{1 \leq i \leq k} (-1)^{i-1} \alpha_i\right| \leq \delta, \quad \text{for all integers } k > 0. \end{cases} \quad (2.1)$$

denoted by $\mathcal{Q}(m,\delta)$ the set of all $\delta$-*deviation* sets, $\zeta(m,\delta) = |\mathcal{Q}(m,\delta)|$.

For a number set $X = \{x_i\}_1^m$, $x_1 > x_2 > \cdots > x_m \geq 0$, and $\alpha \in \mathcal{Q}(m,\delta)$, $\alpha = \{\alpha_i\}$, let $c_0 = 0$, $c_k = \sum_{1 \leq i \leq k} \alpha_i$, for $k > 0$, we define

$$\Delta(X,\alpha) = \left|\sum_{k \geq 0} (-1)^k \sum_{c_k < i \leq c_{k+1}} x_i\right|,$$

**Lemma 1.**

$$\Delta(X,\alpha) \leq \delta \cdot x_1. \tag{2.2}$$

*Proof.* If $\alpha_1 = 1$, let $\alpha'_1 = \alpha_2 - 1, \alpha'_k = \alpha_{k+1},$ for $k > 1$, and $X' = X \setminus \{x_1, x_2\}$, then by the induction

$$\Delta(X,\alpha) \leq (x_1 - x_2) + \Delta(X',\alpha') \leq (x_1 - x_2) + \delta x_3 \leq \delta \cdot x_1.$$

So, assume $\alpha_1 > 1$. Let $\delta_k = \sum_{1 \leq i \leq k} (-1)^{i-1} \alpha_i$. If $\delta_k \geq 0,$ for all $k > 0$, then let $\alpha'_1 = \alpha_1 - 1,$ , $\alpha'_k = \alpha_k$ for all $k > 1$, $X' = X \setminus \{x_1\}$, it is easy to know $\{\alpha'_k\}$ is a $(\delta - 1)$-*deviation* set, hence by induction, it has

$$\Delta(X,\alpha) \leq x_1 + \Delta(X',\alpha') \leq x_1 + (\delta - 1)x_2 \leq \delta \cdot x_1$$

Now suppose that $\delta_s$ is the first one with $\delta_k < 0$, let $\alpha_s = \alpha'_s + \alpha''_s, \alpha'_s = \delta_{s-1}$, $\alpha''_s = (\alpha_s - \delta_{s-1}), \overline{c}_s = c_{s-1} + \delta_{s-1}$, then

$$\Delta(X,\alpha) = |(\sum_{0 \leq k \leq s-1} (-1)^k \sum_{c_k < i \leq c_{k+1}} x_i + (-1)^{s-1} \sum_{1 \leq i \leq \alpha'_s} x_{c_{s-1}+i})$$
$$+ ((-1)^{s-1} \sum_{1 \leq i \leq \alpha''_s} x_{c_s - i} + \sum_{s \leq k} (-1)^k \sum_{c_k < i \leq c_{k+1}} x_i)|$$
$$= |(\sum_{0 \leq k \leq s-1} (-1)^k \sum_{c_k < i \leq c_{k+1}} (x_i - x_{\overline{c}_s}) + (-1)^{s-1} \sum_{1 \leq i \leq \alpha'_s} (x_{c_{s-1}+i} - x_{\overline{c}_s}))$$
$$+ ((-1)^{s-1} \sum_{1 \leq i \leq \alpha''_s} x_{c_s - i} + \sum_{s \leq k} (-1)^k \sum_{c_k < i \leq c_{k+1}} x_i)|$$
$$\leq |(\sum_{0 \leq k \leq s-1} (-1)^k \sum_{c_k < i \leq c_{k+1}} (x_i - x_{\overline{c}_s}) + (-1)^{s-1} \sum_{1 \leq i \leq \alpha'_s} (x_{c_{s-1}+i} - x_{\overline{c}_s}))|$$
$$+ |((-1)^{s-1} \sum_{0 \leq i \leq \alpha''_s} x_{c_s - i} + \sum_{s \leq k} (-1)^k \sum_{c_k < i \leq c_{k+1}} x_i)|$$
$$\leq \delta \cdot (x_1 - x_{\overline{c}_s}) + \delta \cdot x_{\overline{c}_s} = \delta \cdot x_1.$$

$\square$

**Lemma 2.**

$$\zeta(m,1) = 2^{[(m-1)/2]}, \tag{2.3}$$

$$\zeta(m,2) = 3^{[(m-1)/2]}, \tag{2.4}$$

$$\zeta(m,3) = \left((2+\sqrt{2})^{[(m+1)/2]} + (2-\sqrt{2})^{[(m+1)/2]}\right)/4, \tag{2.5}$$

$$\zeta(m,4) = (((5+\sqrt{5})/2)^{[(m+1)/2]} + ((5-\sqrt{5})/2)^{[(m+1)/2]})/5, \tag{2.6}$$

$$\zeta(m,5) = ((2+\sqrt{3})^{[(m+1)/2]} + (2-\sqrt{3})^{[(m+1)/2]} + 2^{[(m+1)/2]})/6. \tag{2.7}$$

*Proof.* Divide $\mathcal{Q}(m,1)$ into two parts

$$\mathcal{Q}_{(1,1)} = \{\alpha \mid \alpha = \{\alpha_i\} \in \mathcal{Q}(m,1), \alpha_1 = 1, \alpha_2 = 1\},$$
$$\mathcal{Q}_{(1,2)} = \{\alpha \mid \alpha = \{\alpha_i\} \in \mathcal{Q}(m,1), \alpha_1 = 1, \alpha_2 = 2\}.$$

Let

$$\mathcal{Q}' = \{\alpha' \mid \alpha = (1,1,\alpha'), \alpha \in \mathcal{Q}_{(1,1)}\},$$
$$\mathcal{Q}'' = \{\alpha'' \mid \alpha = (1,2,\alpha'), \alpha \in \mathcal{Q}_{(1,2)}, \alpha'' = (1,\alpha')\}.$$

It is clear that

$$\mathcal{Q}' = \mathcal{Q}(m-2,\delta), \quad \mathcal{Q}'' = \mathcal{Q}(m-2,\delta)$$

i.e.

$$\zeta(m,1) = 2 \cdot \zeta(m-2,1). \tag{2.8}$$

For (2.4), similarly $\mathcal{Q}(m,2)$ may be divided into three parts

$$\mathcal{Q}(m,2) = \mathcal{Q}_1 \cup \mathcal{Q}_2 \cup \mathcal{Q}_3,$$

where

$$\mathcal{Q}_1 = \{\alpha \mid \alpha = \{\alpha_i\} \in \mathcal{Q}(m,2), \alpha_1 = 1, \alpha_2 = 1\},$$
$$\mathcal{Q}_2 = \{\alpha \mid \alpha = \{\alpha_i\} \in \mathcal{Q}(m,2), \alpha_1 = 1, \alpha_2 = 2, \text{ or } 3\},$$
$$\mathcal{Q}_3 = \{\alpha \mid \alpha = \{\alpha_i\} \in \mathcal{Q}(m,2), \alpha_1 = 2\}.$$

It is easy to know that

$$|\mathcal{Q}_i| = \zeta(m-2,2), \quad i = 1,2,3.$$

So,

$$\zeta(m,2) = 3 \cdot \zeta(m-2,2). \tag{2.9}$$

For (2.5), (2.6) and (2.7) by a similar but a little more investigation, there are following recurrences

$$\zeta(m,3) = 4 \cdot \zeta(m-2,3) - 2 \cdot \zeta(m-4,3), \tag{2.10}$$

$$\zeta(m,4) = 5 \cdot \zeta(m-2,4) - 5 \cdot \zeta(m-4,4), \tag{2.11}$$

$$\zeta(m,5) = 6 \cdot \zeta(m-2,5) - 9 \cdot \zeta(m-4,5) + 2 \cdot \zeta(m-6,5). \tag{2.12}$$

Then formula (2.5), (2,6) and (2.7) are followed by the generating functions with recurrences (2.10), (2.11) and (2.12). □

Suppose that $\mathcal{F}$ is a set of integer numbers contained in $[0,n]$, $\binom{\mathcal{F}}{k}$ as usual stand for the set

of all $k$-subset of $\mathcal{F}$. We define

$$\Omega_k(\mathcal{F},\delta) = \left\{ \Delta(X,\alpha) \mid X \in \binom{\mathcal{F}}{k}, \alpha \in \mathcal{Q}(k,\delta) \right\}.$$

Then,

$$|\Omega_k(\mathcal{F},\delta)| = \binom{|\mathcal{F}|}{k} \cdot |\mathcal{Q}(k,\delta)|. \tag{2.13}$$

Suppose that $\mathfrak{S}$ is a $B_k-sequence$ contained in $[0,n]$, by Lemma 1, each element in $\Omega_k(\mathfrak{S},\delta)$ is a positive integer no greater than $\delta \cdot n$.

Moreover, for two pairs of subsets of $\mathfrak{S}$, $\{A_1,B_1\}$ and $\{A_2,B_2\}$, $A_i \cap B_i = \varnothing$, $|A_i|+|B_i| \leq k$, $i=1,2$, and $|A_1|=|A_2|$, $|B_1|=|B_2|$, if

$$\sum_{x_1 \in A_1} x_1 - \sum_{y_1 \in B_1} y_1 = \sum_{x_2 \in A_2} x_2 - \sum_{y_2 \in B_2} y_2.$$

then $A_1 = A_2$, $B_1 = B_2$. Hence we know the numbers in $\Omega_{2m}(\mathfrak{S},\delta)$ all are different. As to $\Omega_{2m-1}(\mathfrak{S},\delta)$, for $\sum_i \alpha_{2i-1} - \sum_i \alpha_{2i} = 1$, or $-1$, so the numbers in $\Omega_{2m-1}(\mathfrak{S},\delta)$ may be classified into two parts, one contains the numbers with $m$ positive items and $(m-1)$ negative items, the other part contains the numbers with $(m-1)$ positive items and $m$ negative items, clearly, the numbers in each part are different, so at least half of them are different, hence,

$$\begin{array}{l}\binom{|\mathfrak{S}|}{2m} \cdot \zeta(2m,\delta) \leq \delta \cdot n. \\ \binom{|\mathfrak{S}|}{2m-1} \cdot \dfrac{\zeta(2m-1,\delta)}{2} \leq \delta \cdot n.\end{array} \tag{2.14}$$

Besides, for two non-negative integers $s$, $k$, $k \leq s$, as usual denoted by $[s]_k = s \cdot (s-1) \cdots \cdot (s-k+1)$, there is the estimation

**Lemma 3.**

$$[s]_k \geq (s - \mu \cdot k)^k, \qquad \mu \leq (e-1)/e. \tag{2.15}$$

*Proof.* It is easy to demonstrate that $[s]_k / (s-((e-1)k/e))^k$ is monotonic increased when

$s \leq \left(1 + \dfrac{1}{e(e-1)}\right) \cdot k$ and monotonic decreased when $s > \left(1 + \dfrac{1}{e(e-1)}\right) \cdot k$. For

$$\lim_{s \to \infty} \dfrac{[s]_k}{\left(s - ((e-1)k/e)\right)^k} = 1,$$

and Stirling's formula $k! \geq \sqrt{2\pi k}\,(k/e)^k$, hence

$$\dfrac{[s]_k}{\left(s - ((e-1)k/e)\right)^k} \geq \min\left\{1, \dfrac{k!}{(k/e)^k}\right\} \geq 1. \qquad \square$$

Consequently, we obtain

**Theorem 1.**

$$\Phi_k(n) \leq \left(\dfrac{\delta \cdot (1+\tau) \cdot k!}{\zeta(k,\delta)}\right)^{1/k} \cdot n^{1/k} + \mu \cdot k. \qquad (2.16)$$

$\tau = 0$ or $1$, as $k$ even or odd, $\mu \leq \dfrac{e-1}{e}$.

*Note.* In fact, for $s \geq k^2/4$ and for each $i < k/2$, there is

$$(s-i)(s-k+i+1) \geq \left(s - \dfrac{k}{2}\right)^2.$$

It follows,

$$[s]_k \geq \left(s - \dfrac{k}{2}\right)^k, \qquad \text{for } s \geq k^2/4.$$

For the numbers $\zeta(m,\delta)$ play an important role in the result above, we give some further investigation. In general, the recurrence of sequence $\zeta(m,\delta)$ may be written as

$$\zeta(m,\delta) + \sum_{i \geq 1} a_{\delta,i} \cdot \zeta(m-2i,\delta) = 0.$$

Denoted by $p_\delta(x)$ the characteristic polynomials correspond to the recurrence above.

At first we give a general expression for the deduction procedure applied in Lemma 2.

**Lemma 4.** For a integer $k, 0 < k \leq \delta$, let $\mathcal{Q}_k(m,\delta) = \{\alpha \mid \alpha = \{\alpha_i\} \in \mathcal{Q}(m,\delta), \alpha_1 = k\}$, then

$$|\mathcal{Q}_k(m,\delta)| = \begin{cases} 2 \times |\mathcal{Q}(m-2,\delta)|, & \text{for } k=1, \\ \left| \mathcal{Q}(m-2,\delta) \setminus \bigcup_{1 \le i \le k-2} \mathcal{Q}_i(m-2,\delta) \right|, & \text{for } k>1. \end{cases} \quad (2.17)$$

*Proof.* Divide $\mathcal{Q}_k(m,\delta)$ into two parts, $\mathcal{Q}_k(m,\delta) = \mathcal{I} \cup \mathcal{J}$,

$$\mathcal{I} = \{\alpha \mid \alpha = \{\alpha_i\} \in \mathcal{Q}_k(m,\delta), \alpha_2 = 1\}, \quad \mathcal{J} = \{\alpha \mid \alpha = \{\alpha_i\} \in \mathcal{Q}_k(m,\delta), \alpha_2 > 1\}.$$

For set $\mathcal{J}$, let $\mathcal{J}' = \{(\alpha_1-1, \alpha_2-1, \alpha_3, \alpha_4, \cdots) \mid \alpha = \{\alpha_i\} \in \mathcal{J}\}$. Clearly, $|\mathcal{J}| = |\mathcal{J}'|$, and $\mathcal{J}' = \mathcal{Q}_{k-1}(m-2,\delta)$.

For set $\mathcal{I}$, let $\mathcal{I}' = \{(\alpha_3', \alpha_4, \alpha_5, \cdots) \mid \alpha = \{\alpha_i\} \in \mathcal{I}, \alpha_3' = \alpha_1 - \alpha_2 + \alpha_3\}$. For $\alpha_1 = k, \alpha_2 = 1$, and $\alpha_3' \ge k$, so $\mathcal{I}' = \bigcup_{j \ge k} \mathcal{Q}_j(m-2,\delta)$. Besides, it is easy to know $|\mathcal{I}| = |\mathcal{I}'|$. Hence, for $k=1$, $|\mathcal{I}| = |\mathcal{J}| = |\mathcal{Q}(m-2,\delta)|$, and $|\mathcal{Q}_1(m,\delta)| = 2 \times |\mathcal{Q}(m-2,\delta)|$, and for $k>1$,

$$|\mathcal{Q}_k(m,\delta)| = |\mathcal{I} \cup \mathcal{J}| = \left| \bigcup_{j \ge k-1} \mathcal{Q}_j(m-2,\delta) \right| = \left| \mathcal{Q}(m-2,\delta) \setminus \bigcup_{j \le k-2} \mathcal{Q}_j(m-2,\delta) \right|.$$

□

With (2.17) and by the induction, it has

**Corollary 1.** For any positive integers $\delta$ and $k$, $k \le \delta$, there are constants $\{a_{k,i}\}_{i \ge 1}$ such that

$$\left| \bigcup_{j \le k} \mathcal{Q}_j(m,\delta) \right| + \sum_{i \ge 1} a_{k,i} \cdot \zeta(m-2i,\delta) = 0. \quad (2.18)$$

For the characteristic polynomials $p_\delta(x)$ we have

**Theorem 2.**

$$p_\delta(x) = \left( \frac{1+\sqrt{1-4x}}{2} \right)^{\delta+1} + \left( \frac{1-\sqrt{1-4x}}{2} \right)^{\delta+1}. \quad (2.19)$$

*Proof.* Let $\tilde{\mathcal{Q}} = \bigcup_{1 \le i \le \delta-1} \mathcal{Q}_i(m,\delta)$, by Corollary 1,

$$|\tilde{\mathcal{Q}}| + \sum_{i \ge 1} a_{(\delta-1),i} \cdot \zeta(m-2i,\delta) = 0,$$

and

$$|Q(m,\delta)| = \left|\bigcup_{1 \leq i \leq \delta-1} Q_i(m,\delta) \cup Q_\delta(m,\delta)\right| = \left|\tilde{Q} \cup Q_{\delta-1}(m-2,\delta) \cup Q_\delta(m-2,\delta)\right|$$

$$= \left|\tilde{Q} \cup (Q(m-2,\delta) \setminus \bigcup_{1 \leq k \leq (\delta-2)} Q_k(m-2,\delta))\right|.$$

Hence,

$$\zeta(m,\delta) = -\sum_{i \geq 1} a_{(\delta-1),i} \cdot \zeta(m-2i,\delta) + \left(\zeta(m-2,\delta) + \sum_{i \geq 1} a_{(\delta-2),i} \cdot \zeta(m-2-2i,\delta)\right) \quad (2.20)$$

$$= -\sum_{i \geq 1} (a_{(\delta-1),i} - a_{(\delta-2),(i-1)}) \cdot \zeta(m-2i,\delta). \quad (a_{(\delta-2),0} = 1)$$

Correspondingly, for the characteristic polynomials $p_\delta(x)$, there is the recurrence

$$p_\delta(x) = p_{\delta-1}(x) - x \cdot p_{\delta-2}(x). \quad (2.21)$$

Accordingly, by the way of generating function, it follows

$$p_\delta(x) = \left(\frac{1+\sqrt{1-4x}}{2}\right)^{\delta+1} + \left(\frac{1-\sqrt{1-4x}}{2}\right)^{\delta+1}.$$

□

**Theorem 3.**

$$\zeta(m,\delta) = \frac{1}{(\delta+1)} \sum_{1 \leq i \leq d} \left(\frac{1}{\vartheta_i}\right)^{[(m+1)/2]}. \quad (2.22)$$

where $\vartheta_i, i = 1, 2, \cdots, d$, $d = \deg(p_\delta(x))$, are the roots of characteristic polynomial $p_\delta(x)$.

The proof of Theorem above will need following auxiliary results.

**Lemma 5.**

$$\sum_k \binom{r}{2k}\binom{k}{t} = \sum_k \binom{r}{2k-1}\binom{k-1}{t} + \sum_k \binom{r-1}{2k-1}\binom{k-1}{t-1}. \quad (2.23)$$

$$\sum_k \binom{r}{2k-1}\binom{k-1}{t} = \sum_k \binom{r}{2k}\binom{k}{t} + \sum_k \binom{r-1}{2k}\binom{k}{t-1}. \quad (2.24)$$

*Proof.* Let

$$f(x) = \left(1+\sqrt{1+x}\right)^r, \quad g(x) = \left(\sqrt{1+x}-1\right) \cdot f(x).$$

Suppose that

$$f(x) = \sum_k a_k x^k + \sum_k b_k x^k \sqrt{1+x}, \quad g(x) = \sum_k c_k x^k + \sum_k d_k x^k \sqrt{1+x},$$

It is easy to know

$$a_t = \sum_k \binom{r}{2k}\binom{k}{t}, \quad b_t = \sum_k \binom{r}{2k+1}\binom{k}{t}.$$

And,

$$c_t = b_{t-1} - a_t = \sum_k \binom{r}{2k-1}\binom{k-1}{t} - \sum_k \binom{r}{2k}\binom{k}{t}$$

$$d_t = a_t - b_t = \sum_k \binom{r}{2k}\binom{k}{t} - \sum_k \binom{r}{2k+1}\binom{k}{t}$$

On the other hand,

$$g(x) = x \cdot \left(1 + \sqrt{1+x}\right)^{r-1}$$

Hence,

$$c_t = \sum_k \binom{r-1}{2k}\binom{k}{t-1}, \quad d_t = \sum_k \binom{r-1}{2k+1}\binom{k}{t-1}$$

It follows

$$\sum_k \binom{r}{2k}\binom{k}{t} = \sum_k \binom{r}{2k-1}\binom{k-1}{t} + \sum_k \binom{r-1}{2k-1}\binom{k-1}{t-1},$$

$$\sum_k \binom{r}{2k-1}\binom{k-1}{t} = \sum_k \binom{r}{2k}\binom{k}{t} + \sum_k \binom{r-1}{2k}\binom{k}{t-1}.$$

$\square$

**Corollary 2.**

$$\frac{t}{r+1} a_{r,t} = \frac{t}{r} a_{r-1,t} - \frac{t-1}{r-1} a_{r-2,t-1}. \tag{2.25}$$

*Proof.* By Theorem 2, we know

$$a_{r,t} = \frac{(-4)^t}{2^r} \sum_k \binom{r+1}{2k}\binom{k}{t}.$$

So, with use of (2.23),

$$\frac{t}{r} a_{r-1,t} - \frac{t-1}{r-1} a_{r-2,t-1} = \frac{1}{2^{r-1}} \sum_k \frac{t}{r}\binom{r}{2k}\binom{k}{t}(-4)^t - \frac{1}{2^{r-2}} \sum_k \frac{t-1}{r-1}\binom{r-1}{2k}\binom{k}{t-1}(-4)^{t-1}$$

$$= \frac{(-4)^t}{2^{r+1}} \left( \sum_k 2 \times \binom{r-1}{2k-1}\binom{k-1}{t-1} + \sum_k \binom{r-2}{2k-1}\binom{k-1}{t-2} \right)$$

$$= \frac{(-4)^t}{2^{r+1}} \sum_k \binom{r}{2k-1}\binom{k-1}{t-1} - \sum_k \binom{r-1}{2k-2}\binom{k-1}{t-1}$$

$$+ \sum_k \binom{r-1}{2k-1}\binom{k-1}{t-1} + \sum_k \binom{r-2}{2k-1}\binom{k-1}{t-2}$$

$$= \frac{(-4)^t}{2^{r+1}} \sum_k \binom{r}{2k-1}\binom{k-1}{t-1}$$

$$= \frac{t}{r+1} a_{r,t}.$$

□

Write $\zeta(2r,\delta)$ ( or, $\zeta(2r-1,\delta)$ ) as $\omega_r(\delta)$, simply as $\omega_r$ if there is no rise of confusions.

**Corollary 3.** Suppose that $r$ and $\delta$ are two positive integers, denoted by $a_{\delta,0}=1$, $\omega_0 = \dfrac{r}{\delta+1}$, then

$$\sum_{0 \le i \le r} \omega_{r-i}(\delta) \cdot a_{\delta,i} = 0. \tag{2.26}$$

*Proof.* At first we note that $a_{\delta,i}=0$, if $i > d = \deg(p_\delta(x)) = [(\delta+1)/2]$, so for $r > d$, $\sum_{0 \le i \le r} \omega_{r-i}(\delta) \cdot a_{\delta,i} = \sum_{0 \le i \le d} \omega_{r-i}(\delta) \cdot a_{\delta,i} = 0$, for which is just right the recurrence correlation satisfied by the sequence $\{\omega_r(\delta)\}_{r \ge 1}$. Hence, it is assumed $r \le d$.

We take induction on $r$ and $\delta$. Clearly, (2.26) is true for $r=1$, as $a_{\delta,1}=-(\delta+1)$. By (2.3), we know $\omega_r(1) = 2^{r-1}$, so (2.26) hold for $\delta=1$. Besides, it is easy to know that for $r \le \delta$, $\omega_r(\delta)=\omega_r(\delta-1)=\cdots=\omega_r(r)$. Take use of $a_{\delta,r} = a_{\delta-1,r} - a_{\delta-2,r-1}$, and Corollary 2, it has,

$$\sum_{0 \le i \le r} \omega_{r-i} \cdot a_{\delta,i} = \omega_r + \sum_{1 \le i \le r-1} \omega_{r-i} \cdot (a_{\delta-1,i} - a_{\delta-2,i-1}) + \omega_0 \cdot a_{\delta,r}$$

$$= \sum_{0 \le i \le r-1} \omega_{r-i} \cdot a_{\delta-1,i} - \sum_{1 \le i \le r-1} \omega_{r-i} \cdot a_{\delta-2,i-1} + \frac{r}{\delta+1} a_{\delta,r}$$

$$= \left( \sum_{0 \le i \le r-1} \omega_{r-i} \cdot a_{\delta-1,i} + \frac{r}{\delta} a_{\delta-1,r} \right) - \left( \sum_{1 \le i \le r-1} \omega_{r-i} \cdot a_{\delta-2,i-1} + \frac{r-1}{\delta-1} a_{\delta-2,r-1} \right)$$

$$= 0.$$

□

*The Proof of Theorem 3.* We take the induction to prove formula (2.22). For positive integer $k$, denoted by

$$\xi_k = \sum_i \vartheta_i^{-k}.$$

We know $\vartheta_i^{-1}, i=1,2,\cdots,d$, $d = \deg(p_\delta(x))$, are the roots of the polynomial $x^d \cdot p_\delta(x^{-1})$.

With Newton's formula for the two kinds of symmetric polynomials

$$\begin{cases} \sum_{0 \le k \le r-1} \xi_{r-k} a_{\delta,k} + r a_{\delta,r} = 0, & \text{for } r \le d, \\ \sum_{0 \le k \le d} \xi_{r-k} a_{\delta,k} = 0, & \text{for } r > d. \end{cases} \qquad (2.27)$$
$$\qquad (2.28)$$

For $r \le d$, by induction, with (2.27) and Corollary 3, it has

$$\xi_r = -\left( \sum_{1 \le k \le r-1} \xi_{r-k} a_{\delta,k} + r a_{\delta,r} \right) = -\left( \sum_{1 \le k \le r-1} (\delta+1) \omega_{r-k} a_{\delta,k} + r a_{\delta,r} \right)$$

$$= -(\delta+1)\left( \sum_{1 \le k \le r-1} \omega_{r-k} a_{\delta,k} + \frac{r}{\delta+1} a_{\delta,r} \right)$$

$$= (\delta+1)\omega_r$$

i.e.

$$\omega_r = \frac{1}{(\delta+1)} \xi_r.$$

The case $r > d$ may be proved by (2.28) and the induction. □

**Theorem 4.** Suppose that $\vartheta_1 < \vartheta_2 < \cdots < \vartheta_d$, $d = \deg(p_\delta(x))$, are the zeros of the characteristic polynomials $p_\delta(x)$, then

$$\vartheta_k = \frac{1}{4\cos^2((k-0.5)\pi/(\delta+1))}, \qquad k = 1, 2, \cdots, d. \qquad (2.29)$$

*Proof.* From the expression (2.19) of the polynomials $p_\delta(x)$, we can know the zeros of $p_\delta(x)$ all are positive numbers greater than $1/4$, suppose that $\vartheta$ is a root of $p_\delta(x)$, then it may be written as

$$\vartheta = \frac{1+b^2}{4}.$$

$b$ is a positive number. Suppose that $\cos\theta = 1/\sqrt{1+b^2}$, then

$$p_\delta(\vartheta) = \left(\frac{1+b\cdot i}{2}\right)^{\delta+1} + \left(\frac{1-b\cdot i}{2}\right)^{\delta+1} = \left(\frac{\sqrt{1+b^2}}{2}e^{\theta i}\right)^{\delta+1} + \left(\frac{\sqrt{1+b^2}}{2}e^{-\theta i}\right)^{\delta+1}$$

$$= \left(\frac{\sqrt{1+b^2}}{2}\right)^{\delta+1} \cdot 2\cdot \cos(\delta+1)\theta.$$

It follows that

$$\theta = \frac{(k+(1/2))\pi}{(\delta+1)}, \qquad k \text{ is an integer.}$$

and

$$\vartheta = \frac{1}{4\cos^2((k+0.5)\pi/(\delta+1))}.$$

Hence,

$$\vartheta_k = \frac{1}{4\cos^2((k-0.5)\pi/(\delta+1))}, \qquad k=1,2,\cdots,d.$$

$\square$

As an application, we have

**Theorem 5.**

$$\Phi_k(n) \leq \left(\frac{e\cdot \pi^2 \cdot [(k+1)/2]}{2^{(k+2)}} \cdot k!\right)^{1/k} \cdot n^{1/k} + \frac{e-1}{e}k. \qquad (2.30)$$

*Proof.* That (2.30) is true for $k \leq 36$ may be verified directly by Theorem 1 and Lemma 2, so assume $k > 36$.

From the Theorem 4, we know that the smallest root of $p_\delta(x)$, $\vartheta_1 = \frac{1}{4\cos^2(\pi/(2(\delta+1)))}$, and

it is easy to know that

$$\frac{1}{\cos^2 x} \leq 1 + x^2 + x^4, \quad \text{for } |x| \leq \frac{1}{\sqrt{3}},$$

hence,

$$\vartheta_1 \leq \frac{1+\sigma}{4}, \qquad \sigma = (\pi/2(\delta+1))^2 + (\pi/2(\delta+1))^4.$$

By Theorem 1 and 2, it has

$$\zeta(m,\delta) \geq \frac{1}{(\delta+1)}\left(\frac{4}{1+\sigma}\right)^{[(m+1)/2]},$$

and

$$\Phi_k(n) \leq \left( \left( \frac{1+\sigma}{4} \right)^{[(k+1)/2]} \cdot (\delta+1) \cdot \delta \cdot (1+\tau) \cdot k! \right)^{1/k} \cdot n^{1/k} + \frac{e-1}{e} k.$$

where $\tau = 0$ or $1$, as $k$ even or odd.

Take

$$\delta = \left[ \frac{\pi}{2} [(k+3)/2]^{1/2} \right] - 1,$$

It gives,

$$\Phi_k(n) \leq \left( \frac{e \cdot \pi^2 \cdot (1+\tau)}{4^{[(k+3)/2]}} \cdot [(k+1)/2] \cdot k! \right)^{1/k} \cdot n^{1/k} + \frac{e-1}{e} k.$$

and the proof of Theorem 5 is completed. $\square$

In the following, we take two examples to compare the estimations (2.16), with (1.11) and (1.9).

For example $k = 60$, (1.9) and (1.11) give respectively

$$\Phi_{60}(n) \leq \left( 2.11 \times 10^{66} \right)^{1/60} \cdot n^{1/60} + o(n^{1/60}),$$

$$\Phi_{60}(n) \leq \left( 1.227 \times 10^{66} \right)^{1/60} \cdot n^{1/60} + o(n^{1/60}).$$

In (2.16) take $\delta = 7$, it gives

$$\Phi_{60}(n) \leq \left( 1.29 \times 10^{66} \right)^{1/60} \cdot n^{1/60} + 38.$$

For another example $k = 300$, (1.9) and (1.11) are respectively

$$\Phi_{300}(n) \leq \left( 4.896 \times 10^{527} \right)^{1/300} \cdot n^{1/300} + o(n^{1/300}),$$

$$\Phi_{300}(n) \leq \left( 1.315 \times 10^{527} \right)^{1/300} \cdot n^{1/300} + o(n^{1/300}).$$

Take $\delta = 18$, (2.16) gives

$$\Phi_{300}(n) \leq \left( 1.435 \times 10^{527} \right)^{1/300} \cdot n^{1/300} + 190.$$

Next, we will provide some data of test to verify the results of Theorem 2, 3 and 4.

Suppose that $\alpha$ is a odd positive integer less than $2^m$, write it in binary, and add some bits "0" to the end to be a array of of $m$ bits if it less than $m$ bits. We call the adjacent bits with same digits as a block, denoted by $\alpha_i$ the number of bits in the $i$-th block, $i = 1, 2, \cdots$. For $\alpha$ is odd, the first block should be the one with digits "1" and the second block should be the one with digits "0", and so on, which is like as $\underbrace{\underbrace{1\cdots1}_{\alpha_1}\underbrace{0\cdots0}_{\alpha_2}\underbrace{1\cdots1}_{\alpha_3}\cdots\cdots}^{m}$. In this way, a member $\alpha \in \mathcal{Q}(m,\delta)$ may be represented as a array of $m$ bits "0" or "1", i.e. a odd integer less than $2^m$. Obviously, the representation is injective. Therefore, to find $\zeta(m,\delta)$ is suffice to find the number of odd integers less than $2^m$ satisfying equation (2.1). This may be implemented easily by a computer with a short programming when $m$ and $\delta$ are not much great. We have checked for $m$ up to 30, $\delta$ from 1 to 12, and calculated $\zeta(m,\delta)$ with the formulas in Theorem 2,3 and 4 to verify the correctness of the results, which are denoted as $\sim\zeta(m,\delta)$ for temporary in order to distinguish the ones above by counting. The test report demonstrates that the results are correct. A part of data is listed in the following

| $m$ | $\delta$ | $\zeta(m,\delta)$ | $\sim\zeta(m,\delta)$ |
|---|---|---|---|
| 30 | 1 | 16384 | 16384.000000 |
| 30 | 2 | 4782969 | 4782968.999997 |
| 30 | 3 | 24963200 | 24963199.999992 |
| 30 | 4 | 47656250 | 47656249.999990 |
| 30 | 5 | 63255670 | 63255669.999991 |
| 30 | 6 | 71705865 | 71705864.999993 |
| 30 | 7 | 75522960 | 75522959.999994 |
| 30 | 8 | 76964985 | 76964984.999994 |
| 30 | 9 | 77416254 | 77416253.999995 |
| 30 | 10 | 77531355 | 77531354.999995 |
| 30 | 11 | 77554700 | 77554699.999995 |
| 30 | 12 | 77558325 | 77558324.999995 |

List 1

The reference code of the program has been attached in the end of the paper as an appendix.

By the way, we present a by-result.

Suppose that $\alpha = \{\alpha_i\} \in \mathcal{Q}(m,\delta)$, for any positive integer $k$, clearly it has

$$\left|\sum_{1\leq i\leq k}(-1)^{i-1}\alpha_i\right| \leq \max\left\{\left|\sum_{1\leq i\leq \lceil k/2\rceil}\alpha_{2i-1}\right|,\left|\sum_{1\leq i\leq [k/2]}\alpha_{2i}\right|\right\} \leq \lceil m/2\rceil.$$

It follows that $\alpha \in \mathcal{Q}(2r,\delta) \Rightarrow \alpha \in \mathcal{Q}(2r,r)$ if $\delta \geq r$, and so $\zeta(2r,r) = \zeta(2r,\delta)$ if $\delta \geq r$. By Theorem 2,3 and 4, we obtain the following identity,

**Theorem 6.** Suppose that $r, \lambda$ and $\delta$ are three natural numbers, $r \leq \delta \leq \lambda$, then there is

$$\frac{1}{(\delta+1)} \sum_{1\le i\le \lceil \delta/2\rceil} \cos^{2r}\big((2i-1)\pi/(2\delta+2)\big) = \frac{1}{(\lambda+1)} \sum_{1\le i\le \lceil \lambda/2\rceil} \cos^{2r}\big((2i-1)\pi/(2\lambda+2)\big).$$

(2.31)

let $\lambda \to \infty$ in the right-hand side of equation (2.31), then it follows

$$\frac{1}{(\delta+1)} \sum_{1\le i\le \lceil \delta/2\rceil} \cos^{2r}\big((2i-1)\pi/(2\delta+2)\big) = \int_0^{1/2} \cos^{2r} x\pi dx = \frac{1}{2^{2r+1}} C_{2r}^r,$$

(2.32)

and so,

$$\zeta(2r,r) = \frac{1}{2} C_{2r}^r.$$

(2.33)

In fact, the result (2.33) may be obtained more simply by the integer representation of the members in $\mathcal{Q}(m,\delta)$ mentioned above, it is clear that $\zeta(2r,r)$ is equal to the number of odd integers less than $2^{2r}$ with $r$ bits "1", hence

$$\zeta(2r,r) = C_{2r-1}^{r-1} = \frac{1}{2} C_{2r}^r.$$

In the following, we will present another formula for $\zeta(2r,\delta)$.

**Theorem 7.**

$$\zeta(2r,\delta) = \frac{1}{2} C_{2r}^r + \sum_{i>0}(-1)^i \cdot C_{2r}^{r-i\cdot(\delta+1)}.$$

(2.34)

In the proof of Theorem 7, it will be applied the following formulas

a) $\cos^{2r} x = \dfrac{1}{2^{2r}} C_{2r}^r + \dfrac{1}{2^{2r-1}} \sum_{1\le k\le r} C_{2r}^{r-k} \cos 2kx.$ (2-a)

b) $\sum_{1\le i\le \lambda} \cos\big((2i-1)k\pi/(2\lambda)\big) = \begin{cases} (-1)^{k/(2\lambda)} \lambda, & \text{if } (2\lambda)\mid k, \\ 0, & \text{otherwise.} \end{cases}$ (2-b)

c) $\sum_{1\le i\le \lambda} \cos\big((2i-1)k\pi/(2\lambda+1)\big) = \begin{cases} (-1)^{k-1}\dfrac{1}{2}, & \text{if } (2\lambda+1)\nmid k, \\ (-1)^{k/(2\lambda+1)} \lambda, & \text{if } (2\lambda+1)\mid k. \end{cases}$ (2-c)

d) $\dfrac{1}{2} C_{2r}^r + \sum_{1\le k\le r}(-1)^k \cdot C_{2r}^{r-k} = 0.$ (2-d)

*The Proof of Theorem 7*:

At first assume that $\delta$ is odd, $\delta = 2\lambda-1,$ then

$$\zeta(2r,\delta) = \frac{4^r}{(\delta+1)} \sum_{1\leq i\leq \lceil \delta/2\rceil} \cos^{2r}\left((2i-1)\pi/(2\delta+2)\right)$$

$$= \frac{4^r}{(2\lambda)} \sum_{1\leq i\leq \lambda} \left( \frac{1}{2^{2r}} C_{2r}^r + \frac{1}{2^{2r-1}} \sum_{1\leq k\leq r} C_{2r}^{r-k} \cos k(2i-1)\pi/(2\lambda) \right)$$

$$= \frac{1}{2} C_{2r}^r + \frac{1}{\lambda} \sum_{1\leq i\leq \lambda} \sum_{1\leq k\leq r} C_{2r}^{r-k} \cos k(2i-1)\pi/(\delta+1)$$

$$= \frac{1}{2} C_{2r}^r + \frac{1}{\lambda} \sum_{1\leq k\leq r} C_{2r}^{r-k} \sum_{1\leq i\leq \lambda} \cos k(2i-1)\pi/(\delta+1)$$

$$= \frac{1}{2} C_{2r}^r + \sum_{i>0} (-1)^i \cdot C_{2r}^{r-i\cdot(\delta+1)}$$

If $\delta$ is even, $\delta = 2\lambda$, then

$$\zeta(2r,\delta) = \frac{4^r}{(\delta+1)} \sum_{1\leq i\leq \lceil \delta/2\rceil} \cos^{2r}\left((2i-1)\pi/(2\delta+2)\right)$$

$$= \frac{4^r}{(2\lambda+1)} \sum_{1\leq i\leq \lambda} \left( \frac{1}{2^{2r}} C_{2r}^r + \frac{1}{2^{2r-1}} \sum_{1\leq k\leq r} C_{2r}^{r-k} \cos k(2i-1)\pi/(\delta+1) \right)$$

$$= \frac{\lambda}{(2\lambda+1)} C_{2r}^r + \frac{2}{(2\lambda+1)} \sum_{1\leq i\leq \lambda} \sum_{1\leq k\leq r} C_{2r}^{r-k} \cos k(2i-1)\pi/(2\lambda+1)$$

$$= \frac{\lambda}{(2\lambda+1)} C_{2r}^r + \frac{2}{(2\lambda+1)} \sum_{1\leq k\leq r} C_{2r}^{r-k} \sum_{1\leq i\leq \lambda} \cos k(2i-1)\pi/(2\lambda+1)$$

$$= \frac{\lambda}{(2\lambda+1)} C_{2r}^r + \frac{2\lambda}{(2\lambda+1)} \sum_{i>0} (-1)^i \cdot C_{2r}^{r-i\cdot(\delta+1)} - \frac{1}{(2\lambda+1)} \sum_{k\nmid (2\lambda+1)} (-1)^k \cdot C_{2r}^{r-k}$$

$$= \frac{1}{2} C_{2r}^r + \sum_{i>0} (-1)^i \cdot C_{2r}^{r-i\cdot(\delta+1)} - \frac{1}{(2\lambda+1)} \cdot \left( \frac{C_{2r}^r}{2} + \sum_i (-1)^i C_{2r}^{r-i\cdot(\delta+1)} + \sum_{k\nmid (2\lambda+1)} (-1)^k C_{2r}^{r-k} \right)$$

$$= \frac{1}{2} C_{2r}^r + \sum_{i>0} (-1)^i \cdot C_{2r}^{r-i\cdot(\delta+1)} - \frac{1}{(2\lambda+1)} \left( \frac{1}{2} C_{2r}^r + \sum_{1\leq k\leq r} (-1)^k C_{2r}^{r-k} \right)$$

$$= \frac{1}{2} C_{2r}^r + \sum_{i>0} (-1)^i \cdot C_{2r}^{r-i\cdot(\delta+1)}.$$

□

With the integer representation of $\delta$-deviation set, we realize that there may be an alternative definition for $\zeta(m,\delta)$.

Suppose that $\alpha$ is a array of $m$ bits of "0" or "1", denoted by $|\alpha|, \langle\alpha\rangle_0$ and $\langle\alpha\rangle_1$ the length of $\alpha$, the number of bits "0" and "1" of array $\alpha$ respectively, and denoted by $\alpha_k$ the segment consist of the first $k$ bits of $\alpha$, for the positive integers $m$ and $\delta$, define $z(m,\delta)$ is the number of the arrays $\alpha$ of $m$ bits meet the following conditions

$$\begin{cases} \alpha_1 = \text{"1"}, \ |\alpha| = m, \\ |\langle\alpha\rangle_0 - \langle\alpha\rangle_1| \leq 1, \ |\langle\alpha_k\rangle_0 - \langle\alpha_k\rangle_1| \leq \delta, \quad \text{for all integer } k > 0. \end{cases} \quad (2.1\text{-a})$$

It is not difficult to know that $z(m,\delta) = \zeta(m,\delta)$, that is, both are equivalent.

With the definition of $\zeta(m,\delta)$ above, for the formula (2.34) there is appeared combinatorial meaning of include-exclude form, so, it is appealing to give an exposition for (2.34) in this way. In the following, we try to give a proof/interpretation for Theorem 7 by include-exclude means.

At first, we introduce some notations will be used in the next.

Denoted by $V = \{0,1\}$, $V^m$ is as usual the set of all arrays of length $m$ over domain $V$. As before, for a $\alpha \in V^m$, denoted by $\alpha_k$ the segment of the first $k$ bits of $\alpha$, $\langle\alpha\rangle_0$ and $\langle\alpha\rangle_1$ the number of bits "0" and "1" of array $\alpha$ respectively. For a $\alpha \in V^m$ and a positive integer $k \leq m$, define $d_k(\alpha) = |\langle\alpha_k\rangle_0 - \langle\alpha_k\rangle_1|$, simply write $d_m(\alpha) = d(\alpha)$.

**Definition 2.** For a positive integer $\delta \leq m$, let $\Im(m,\delta)$ be the set of arrays $\alpha \in V^m$ meet the condition

$$d(\alpha) \leq 1, \quad d_k(\alpha) \leq \delta, \quad \text{for all integer } k > 0. \quad (2.1\text{-b})$$

It is clear that $|\Im(m,\delta)| = 2 \cdot \zeta(m,\delta)$. Simply write $\lambda = \delta + 1$, for a $\alpha \in V^m$, define

$$k_{\max}(\alpha) = \max\{k \mid d_k(\alpha) \equiv \lambda \mod(2\lambda)\}.$$

Suppose that $l$, $t$ are two non-negative integers, $t \leq l \leq m$, define

$$\Gamma(t,l) = \{\alpha \mid \alpha \in V^l, \langle\alpha\rangle_1 = t\}, \text{ and}$$

$$\Lambda(t,l) = \{\alpha \mid \alpha \in \Gamma(t,l), 0 < d_k(\alpha) < 2\lambda, \text{ for all integer } k > 0\}.$$

The arrays in $\Lambda(t,l)$ are sometimes called $\Lambda$-array.

Let $\Delta_1(r,m,\delta) = \Gamma(r,m) \setminus \Im(m,\delta)$. In the following, it is assumed that $m = 2r$.

**Lemma 6.** If $\alpha \in \Delta_1(r,m,\delta)$, then there is a positive integer $k$, $1 \leq k \leq m$, such that $d_k(\alpha) = \lambda$.

*Proof.* From the definition of $\Delta_1(r,m,\delta)$, there is at least a positive integer $k$, $1 \leq k \leq m$, such that $d_k(\alpha) \geq \lambda$. Assume that $k$ is of the greatest one, with the "continuality" of the function $d_k(\alpha)$, and $d_m(\alpha) \leq 1$, so $d_k(\alpha) = \lambda$. □

Here we say that function $d_k(\alpha)$ is continuous means that $|d_{k+1}(\alpha) - d_k(\alpha)| \leq 1$ for all positive integer $k$ and all $\alpha \in V^m$.

For a positive integer $k$, $1 \leq k \leq m$, denoted by $E_k = \{\alpha \mid \alpha \in \Delta_1(r,m,\delta), k_{\max}(\alpha) = k\}$. Divide $E_k = E_k^{(0)} \cup E_k^{(1)}$,

$$E_k^{(0)} = \{\alpha \mid \alpha \in E_k, d_k(\alpha) = \langle \alpha_k \rangle_0 - \langle \alpha_k \rangle_1\}, \quad E_k^{(1)} = \{\alpha \mid \alpha \in E_k, d_k(\alpha) = \langle \alpha_k \rangle_1 - \langle \alpha_k \rangle_0\}.$$

Let $\Delta_1^{(i)}(r,m,\delta) = \bigcup_k E_k^{(i)}$, $i = 0,1$, then $\Delta_1(r,m,\delta) = \Delta_1^{(0)}(r,m,\delta) \cup \Delta_1^{(1)}(r,m,\delta)$. For the symmetry, it has $|\Delta_1^{(0)}(r,m,\delta)| = |\Delta_1^{(1)}(r,m,\delta)|$.

For a $\alpha \in E_k^{(0)}$, then $k = \lambda + 2v$, $\langle \alpha_k \rangle_0 = v + \lambda$, $\langle \alpha_k \rangle_1 = v$ for some integer $v \geq 0$. Let $\alpha_k$ and $\alpha_k'$ be the segments of the first $k$ bits and the last $(m-k)$ bits respectively, i.e. $\alpha = \alpha_k \alpha_k'$, and $\alpha_k' \in \Lambda(r-v, m-k)$, hence it has

$$\Delta_1^{(0)}(r, 2r, \delta) = \bigcup_{0 \leq v \leq r-\lambda} (\Gamma(v, 2v+\lambda) \times \Lambda(r-v, 2r-2v-\lambda)).$$

For a $\alpha \in V^l$, denoted by $\bar{\alpha}$ the complement of $\alpha$, i.e. $\bar{\alpha} \in V^l$, with that $\alpha \wedge \bar{\alpha} = I$ ($= 1\cdots 1$), where operation $\wedge$ is the addition $\mathrm{mod}(2)$ in bits. Denoted by $\bar{\Lambda}(t,l) = \{\bar{\alpha} \mid \alpha \in \Lambda(t,l)\}$. Clearly, $\bar{\Lambda}(t,l) = \Lambda(l-t,l)$, and $|\Lambda(t,l)| = |\Lambda(l-t,l)|$. Let

$$\tilde{\Delta}_1^{(0)}(r, 2r, \delta) = \bigcup_{0 \leq v \leq r-\lambda} (\Gamma(v, 2v+\lambda) \times \Lambda(r-v-\lambda, 2r-2v-\lambda)).$$

There is a map $\phi: \Delta_1^{(0)}(r,2r,\delta) \to \tilde{\Delta}_1^{(0)}(r,2r,\delta)$, $\alpha \mapsto \tilde{\alpha}$, $\alpha = \alpha_k \alpha_k'$, $\tilde{\alpha} = \alpha_k \bar{\alpha}_k'$, $k = k_{\max}(\alpha)$.

**Lemma 7.** Denoted by $D_v = \Gamma(v, 2v+\lambda) \times \Lambda(r-v-\lambda, m-2v-\lambda)$, then $D_i \cap D_j = \varnothing$,

for $i \neq j$.

*Proof.* Suppose that $\alpha \in D_i \cap D_j, \alpha \neq \varnothing$, then there are

$\alpha = \alpha_k \alpha'_k, \alpha_k \in \Gamma(i, 2i + \lambda), k = \lambda + 2i, \alpha'_k \in \Lambda(r - i - \lambda, m - 2i - \lambda)$, and

$\alpha = \alpha_s \alpha'_s, \alpha_s \in \Gamma(j, 2j + \lambda), s = \lambda + 2j, \alpha'_s \in \Lambda(r - j - \lambda, m - 2j - \lambda)$.

Assume $i < j$, then $\alpha_s = \alpha_k \hat{\alpha}$, $\hat{\alpha}$ is the segment of the first $2(j-i)$ bits of $\alpha'_k$, and $\left| \langle \hat{\alpha} \rangle_0 - \langle \hat{\alpha} \rangle_1 \right| = 0$, which is contradict with the fact that $\alpha'_k$ is a $\Lambda$- array. □

Lemma 7 indicates the map $\phi(\alpha)$ is bijective. For $\tilde{\Delta}_1^{(0)}(r, 2r, \delta) \subseteq \Gamma(r - \lambda, 2r)$, hence

$\left| \Delta_1^{(0)}(r, 2r, \delta) \right| = \left| \tilde{\Delta}_1^{(0)}(r, 2r, \delta) \right| \leq C_{2r}^{r-\lambda}$.

In general, for arbitrary integers $i > 0, v \geq 0$, let

$A_v^{(i)} = \Gamma(v, 2v + (2i - 1)\lambda) \times \Lambda(r - (i-1)\lambda - v, m - 2v - (2i - 1)\lambda)$,

$B_v^{(i)} = \Gamma(v, 2v + (2i - 1)\lambda) \times \Lambda(r - i\lambda - v, m - 2v - (2i - 1)\lambda)$,

$\Delta_i^{(0)}(r, m, \delta) = \bigcup_{v \geq 0} A_v^{(i)}$, $\tilde{\Delta}_i^{(0)}(r, m, \delta) = \bigcup_{v \geq 0} B_v^{(i)}$.

$\Delta_i(r, m, \delta)$ will be written as $\Delta_i(r, m)$ simply if parameter $\delta$ is unchanged in the context.

Similar to Lemma 7, there is

**Lemma 8.**

$A_s^{(i)} \cap A_t^{(i)} = \varnothing, \ B_s^{(i)} \cap B_t^{(i)} = \varnothing$, for all positive integers $s \neq t$. (2.35)

Then, we have

**Lemma 9.**

$\left| \Delta_i^{(0)}(r, m) \right| = \left| \tilde{\Delta}_i^{(0)}(r, m) \right|$, for all integers $i \geq 1$. (2.36)

*Proof.*

$\left| \Delta_i^{(0)}(r, m) \right| = \sum_{v \geq 0} \left| A_v^{(i)} \right|$

$= \sum_{v \geq 0} \left| \Gamma(v, 2v + (2i - 1)\lambda) \right| \times \left| \Lambda(r - (i-1)\lambda - v, m - 2v - (2i - 1)\lambda) \right|$

$= \sum_{v \geq 0} \left| \Gamma(v, 2v + (2i - 1)\lambda) \right| \times \left| \Lambda(r - i\lambda - v, m - 2v - (2i - 1)\lambda) \right|$  ($\because \left| \Lambda(t, l) \right| = \left| \Lambda(l - t, l) \right|$)

$$= \sum_{v \geq 0} \left| B_v^{(i)} \right| = \left| \tilde{\Delta}_i^{(0)}(r,m) \right|.$$

□

Besides, there is

**Lemma 10.** For arbitrary integer $i \geq 1$,

$$\Gamma(r - i\lambda, m) = \tilde{\Delta}_i^{(0)}(r,m) \cup \Delta_{i+1}^{(0)}(r,m), \tag{2.37}$$

$$\Delta_{i+1}^{(0)}(r,m) \cap \tilde{\Delta}_i^{(0)}(r,m) = \varnothing. \tag{2.38}$$

*Proof.* For a $\alpha \in \Gamma(r - i\lambda, m)$, then $k_{\max}(\alpha) = 2v + a\lambda$, $\langle \alpha_{k_{\max}} \rangle_1 = v$, $\langle \alpha_{k_{\max}} \rangle_0 = v + a\lambda$, for some integer $v \geq 0$, $a$ is an odd number. For the continuality of function $d_k(\alpha)$ and $d(\alpha) = 2i\lambda$, it infers that $a = 2i - 1$, or $2i + 1$, the possible odd integers nearest $2i$.

If $a = 2i - 1$, i.e.

$$\alpha \in \Gamma(v, 2v + (2i-1)\lambda) \times \Lambda(r - i\lambda - v, 2r - 2v - (2i-1)\lambda) \subseteq \tilde{\Delta}_i^{(0)}(r,m).$$

If $a = 2i + 1$, i.e.

$$\alpha \in \Gamma(v, 2v + (2i+1)\lambda) \times \Lambda(r - i\lambda - v, 2r - 2v - (2i+1)\lambda) \subseteq \Delta_{i+1}^{(0)}(r,m).$$

So (2.37) has been proved. As to (2.38), Suppose that $\alpha \in A_v^{(i+1)} \cap B_q^{(i)}, \alpha \neq \varnothing$, $\alpha$ may be written as $\alpha = \alpha_s \beta'$ and $\alpha = \alpha_t \beta''$, $s = 2v + (2i+1)\lambda$, $t = 2q + (2i-1)\lambda$, $\beta'$, $\beta''$ are $\Lambda$-arrays. If $t \geq s$, then we may write $\alpha_t = \alpha_s \tau$, and

$$d_t(\alpha) = \left| \langle \alpha_s \rangle_0 + \langle \tau \rangle_0 - \langle \alpha_s \rangle_1 - \langle \tau \rangle_1 \right| = \left| (2i+1)\lambda \pm \left( \langle \tau \rangle_0 - \langle \tau \rangle_1 \right) \right| = (2i-1)\lambda.$$

It follows that $0 < \left| \langle \tau \rangle_0 - \langle \tau \rangle_1 \right| \equiv 0 \mod(2\lambda)$, that is contradicted with that $\beta'$ is $\Lambda$-array. The rest proof for the case $s \geq t$ is similar, so omitted. □

Consequently,

$$\zeta(m,\delta) = \frac{1}{2} |\Im(m,\delta)| = \frac{1}{2} \left( |\Gamma(r,2r)| - \left( |\Gamma(r,2r)| - |\Im(m,\delta)| \right) \right)$$

$$= \frac{1}{2} \left( |\Gamma(r,2r)| - \left( |\Delta_1^{(0)}(r,2r)| + |\Delta_1^{(1)}(r,2r)| \right) \right) = \frac{1}{2} \left( |\Gamma(r,2r)| - 2|\Delta_1^{(0)}(r,2r)| \right)$$

$$= \frac{1}{2} \left( |\Gamma(r,2r)| - 2\left( |\Delta_1^{(0)}(r,2r)| + |\Delta_2^{(0)}(r,2r)| \right) + 2\left( |\Delta_2^{(0)}(r,2r)| + |\Delta_3^{(0)}(r,2r)| \right) - \cdots \right)$$

$$= \frac{1}{2} \left( |\Gamma(r,2r)| - 2\left( |\tilde{\Delta}_1^{(0)}(r,2r)| + |\Delta_2^{(0)}(r,2r)| \right) + 2\left( |\tilde{\Delta}_2^{(0)}(r,2r)| + |\Delta_3^{(0)}(r,2r)| \right) - \cdots \right)$$

$$= \frac{1}{2}\left(|\Gamma(r,2r)| + 2\sum_{i>0}(-1)^i \cdot |\Gamma(r-i\lambda, 2r)|\right)$$

$$= \frac{1}{2}C_{2r}^r + \sum_{i>0}(-1)^i \cdot C_{2r}^{r-i\lambda}.$$

In the next, we develop the investigation above to the more general cases.

**Definition 3.** For the non-negative integers $0 \leq \varepsilon \leq \delta \leq m$, let $\mathfrak{D}(m,\varepsilon,\delta)$ and $\mathfrak{S}(m,\delta)$ be the sets of arrays $\alpha \in V^m$ satisfying conditions (2.1-c) and (2.1-d) respectively,

$$d(\alpha) = \varepsilon, \quad d_k(\alpha) \leq \delta, \quad \text{for all integer } k > 0. \tag{2.1-c}$$
$$d_k(\alpha) \leq \delta, \quad \text{for all integer } k > 0. \tag{2.1-d}$$

Denoted by $\theta(x,y)$ the characteristic function of the parity of integers $x$ and $y$, i.e. $\theta(x,y) = 0$, or $1$, as $x$ and $y$ are or not of same parity, i.e. $x + y \equiv \theta(x,y) \mod(2)$.

**Theorem 8.** Suppose that $m = r+s$, $\varepsilon = |r-s|$, $\lambda = \delta + 1$, then

$$|\mathfrak{D}(m,\varepsilon,\delta)| = (1+\text{sgn}(\varepsilon)) \cdot \left(C_m^r + \sum_{i>0}(-1)^i \left(C_m^{r-i\lambda} + C_m^{s-i\lambda}\right)\right). \tag{2.39}$$

Moreover, denoted by $h = [(m-\delta)/2]$, $\theta = \theta(m,\delta)$, then

$$|\mathfrak{S}(m,\delta)| = \sum_{k=\theta}^{\delta} C_m^{h+k} + 2 \cdot \sum_{i>0}(-1)^i \sum_{k=\theta}^{\delta} C_m^{h-i\lambda+k}. \tag{2.40}$$

Let $\Delta_i^{(0)}(r,m)$ be constructed as before. For them here there are analogous results as Lemma 8,9 and 10, but Lemma 9 with a little difference as that

**Lemma 9'.**
$$\left|\Delta_i^{(0)}(r,m)\right| = \left|\tilde{\Delta}_i^{(0)}(s,m)\right|, \quad \text{for all integers } i \geq 1. \tag{2.36'}$$

*Proof of Theorem 8.* By Theorem 7, we know that (2.39) is true for $\varepsilon = 0$, so assumed $r \neq s$. Let $\Gamma = \Gamma(r,m) \cup \Gamma(s,m)$, and $\Delta_1 = \Gamma \setminus \mathfrak{D}(m,\varepsilon,\delta)$, similar to $\Delta_1(r,m)$ in the proof above, divide $\Delta_1 = \Delta_1^{(0)} \cup \Delta_{(1)}^{(1)}$, for the symmetry, it has $\left|\Delta_1^{(0)}\right| = \left|\Delta_1^{(1)}\right|$. Moreover, it is clear that

$\Delta_1^{(0)} = \Delta_1^{(0)}(r,m) \cup \Delta_1^{(0)}(s,m)$. Thus

$$|\mathfrak{D}(m,\varepsilon,\delta)| = (|\Gamma| - (|\Gamma| - |\mathfrak{D}(m,\varepsilon,\delta)|)) = (|\Gamma| - |\Delta_1|)$$

$$= (|\Gamma| - (|\Delta_1^{(0)}| + |\Delta_1^{(1)}|)) = (|\Gamma| - 2 \cdot |\Delta_1^{(0)}|) = (|\Gamma| - 2(|\Delta_1^{(0)}(r,m)| + |\Delta_1^{(0)}(s,m)|))$$

$$= (|\Gamma(r,m)| - 2(|\Delta_1^{(0)}(r,m)| + |\Delta_2^{(0)}(s,m)|) + 2(|\Delta_2^{(0)}(s,m)| + |\Delta_3^{(0)}(r,m)|) - \cdots)$$

$$+ (|\Gamma(s,m)| - 2(|\Delta_1^{(0)}(s,m)| + |\Delta_2^{(0)}(r,m)|) + 2(|\Delta_2^{(0)}(r,m)| + |\Delta_3^{(0)}(s,m)|) - \cdots)$$

$$= (|\Gamma(r,m)| - 2(|\tilde{\Delta}_1^{(0)}(s,m)| + |\Delta_2^{(0)}(s,m)|) + 2(|\tilde{\Delta}_2^{(0)}(r,m)| + |\Delta_3^{(0)}(r,m)|) - \cdots)$$

$$+ (|\Gamma(s,m)| - 2(|\tilde{\Delta}_1^{(0)}(r,m)| + |\Delta_2^{(0)}(r,m)|) + 2(|\tilde{\Delta}_2^{(0)}(s,m)| + |\Delta_3^{(0)}(s,m)|) - \cdots)$$

$$= \left( |\Gamma(r,m)| + |\Gamma(s,m)| + 2\sum_{i>0}(-1)^i (|\Gamma(r-i\lambda,m)| + |\Gamma(s-i\lambda,m)|) \right)$$

$$= C_m^r + C_m^s + 2 \cdot \sum_{i>0}(-1)^i \left( C_m^{r-i\lambda} + C_m^{s-i\lambda} \right).$$

This has proved (2.39). For (2.40), if $m$ is even, then by (2.39), it has

$$|\mathfrak{S}(m,\delta)| = \sum_{0 \leq \varepsilon \leq \delta} |\mathfrak{D}(m,\varepsilon,\delta)|$$

$$= C_m^{(m/2)} + 2 \cdot \sum_{i>0}(-1)^i C_m^{(m/2)-i\lambda} + \sum_{k=1}^{[\delta/2]} \left( C_m^{(m/2)-k} + C_m^{(m/2)+k} \right)$$

$$+ 2 \cdot \sum_{i>0}(-1)^i \sum_{k=1}^{[\delta/2]} \left( C_m^{[(m/2)]-i\lambda-k} + C_m^{[(m/2)]-i\lambda+k} \right)$$

$$= \sum_{k=0}^{2 \cdot [\delta/2]} C_m^{(m/2)-[\delta/2]+k} + 2 \cdot \sum_{i>0}(-1)^i \sum_{k=0}^{2 \cdot [\delta/2]} C_m^{(m/2)-[\delta/2]-i\lambda+k}$$

$$= \sum_{k=0}^{2 \cdot [\delta/2]} C_m^{(m/2)-(\delta-\theta)/2+k} + 2 \cdot \sum_{i>0}(-1)^i \sum_{k=\theta}^{2 \cdot [\delta/2]} C_m^{(m/2)-(\delta-\theta)/2-i\lambda+k}$$

$$= \sum_{k=0}^{2 \cdot [\delta/2]} C_m^{(m-\delta-\theta)/2+k+\theta} + 2 \cdot \sum_{i>0}(-1)^i \sum_{k=\theta}^{2 \cdot [\delta/2]} C_m^{(m-\delta-\theta)/2-i\lambda+k+\theta}$$

$$= \sum_{k=\theta}^{\delta} C_m^{h+k} + 2 \cdot \sum_{i>0}(-1)^i \sum_{k=\theta}^{\delta} C_m^{h-i\lambda+k}.$$

So, (2.40) has been proved in the case $m$ is even, the case for $m$ is odd may be verified in a similar way, which is omitted.

□

The proof above is just a brief one, in which the similar investigation with the above are omitted.

Now, it is natural to ask the possibility with the methods of generating functions as in the Theorems 2, 3 and 4 to calculus $|\mathfrak{D}(m,\varepsilon,\delta)|$ and $|\mathfrak{S}(m,\delta)|$. In fact, it is easy to know that the investigations before are available for $|\mathfrak{D}(m,\varepsilon,\delta)|$ yet, that is, the sequences $\{|\mathfrak{D}(m,\varepsilon,\delta)|\}$

have same characteristic polynomials $p_\delta(x)$ as $\zeta(m,\delta)$, in fact, $\zeta(m,\delta)$ is just a special case $\varepsilon = 0$ of $|\mathfrak{D}(m,\varepsilon,\delta)|$, hence, it has the general formula

$$|\mathfrak{D}(m,\varepsilon,\delta)| = \sum_{1 \le i \le d} b_i(\varepsilon) \cdot \vartheta_i^{-(m-\varepsilon+2)/2}. \tag{2.41}$$

In the discrete cases, the coefficients $b_i(\varepsilon)$'s may be determined by the initial $d$ values of $|\mathfrak{D}(m,\varepsilon,\delta)|$, $m = \varepsilon, \varepsilon+2, \cdots, \varepsilon+2(d-1)$. For instances,

$$|\mathfrak{D}(m,2,3)| = \frac{1}{\sqrt{2}} \cdot \left( \left(2+\sqrt{2}\right)^{m/2} - \left(2-\sqrt{2}\right)^{m/2} \right), \qquad (m \text{ even}).$$

$$|\mathfrak{D}(m,2,4)| = \frac{1+\sqrt{5}}{5} \cdot \left( \frac{5+\sqrt{5}}{2} \right)^{m/2} + \frac{1-\sqrt{5}}{5} \cdot \left( \frac{5-\sqrt{5}}{2} \right)^{m/2}, \qquad (m \text{ even}).$$

We have seen that $b_i$'s are not again same simple as the ones in the case of $\zeta(m,\delta)$. It seems not much easy to find a simple way to know these coefficients in general. Nevertheless, the result of Theorem 8 and the investigation above reveal some message for these coefficients.

**Theorem 9.** Denoted by $\lambda = \delta+1, d = [(\delta+1)/2]$, then

$$|\mathfrak{D}(m,\varepsilon,\delta)| = \frac{2^{m+1+\mathrm{sgn}(\varepsilon)}}{\lambda} \cdot \sum_{1 \le i \le d} \cos\left(\frac{\varepsilon(2i-1)\pi}{2\lambda}\right) \cdot \cos^m\left(\frac{(2i-1)\pi}{2\lambda}\right). \tag{2.42}$$

Let $\theta(m,\delta)$ be the characteristic function of parity of $m$ and $\delta$ as defined in Theorem 8.

**Theorem 10.** Denoted by $\lambda = \delta+1, d = [(\delta+1)/2], \theta = \theta(m,\delta)$, then

$$|\mathfrak{S}(m,\delta)| = \frac{2^{m+1}}{\lambda} \cdot \sum_{1 \le i \le d} (-1)^{i-1} \frac{\cos^{m+\theta}((2i-1)\pi/(2\lambda))}{\sin((2i-1)\pi/(2\lambda))}. \tag{2.43}$$

In the proofs of Theorems 9, 10, it will be used following formulas

e) $\cos^{2t+1} x = \dfrac{1}{2^{2t}} \sum_{0 \le k \le t} C_{2t+1}^{t-k} \cos(2k+1)x.$ \hfill (2-e)

f) $\sum_{1 \le j \le d} \cos\left((2i-1)(2j-1)\pi/(4d)\right) = (-1)^{i-1}/\sin((2i-1)\pi/(4d)).$ \hfill (2-f)

g) $\sum_{1 \le j \le d} \cos\left((2i-1)(2j-1)\pi/(4d+2)\right) = (-1)^{i-1}\cot((2i-1)\pi/(4d+2)).$ \hfill (2-g)

h) $\sum_{0\le j<d}\cos((2i-1)2j\pi/(4d))=\frac{1}{2}(1+(-1)^{i-1}\cot((2i-1)\pi/(4d)))$.   (2-h)

i) $\sum_{0\le j\le d}\cos((2i-1)2j\pi/(4d+2))=\frac{1}{2}(1+(-1)^{i-1}/\sin((2i-1)\pi/(4d+2)))$.   (2-i)

*Proof of Theorem 9:*   Suppose that $m$ is odd, $m=2t+1$, for $\varepsilon$ and $m$ have to be same parity, so assume $\varepsilon=2\varepsilon_0+1$. By formula (2-e), it has

$$\cos(\varepsilon x)\cdot\cos^{2t+1}x=\frac{1}{2^{2t}}\sum_{0\le k\le t}C_{2t+1}^{t-k}\cdot\cos(\varepsilon x)\cdot\cos(2k+1)x$$

$$=\frac{1}{2^{2t+1}}\sum_{0\le k\le t}C_{2t+1}^{t-k}\cdot\left(\cos((2k+1+\varepsilon)x)+\cos((2k+1-\varepsilon)x)\right).$$

If $\delta=2d-1$, then

$$\sum_{1\le i\le d}\cos(\varepsilon((2i-1)\pi/(2\delta+2))\cdot\cos^{2t+1}((2i-1)\pi/(2\delta+2))$$

$$=\frac{1}{2^{2t+1}}\cdot\sum_{0\le k\le t}C_{2t+1}^{t-k}\cdot\sum_{1\le i\le d}\left(\cos\left(\frac{(2i-1)\pi}{(2\delta+2)}((2k+1)+\varepsilon)\right)+\cos\left(\frac{(2i-1)\pi}{(2\delta+2)}((2k+1)-\varepsilon)\right)\right)$$

$$=\frac{1}{2^{2t+1}}\cdot\sum_{0\le k\le t}C_{2t+1}^{t-k}\cdot\sum_{1\le i\le d}\left(\cos\left(\frac{(2i-1)\pi}{2d}(k+1+\varepsilon_0)\right)+\cos\left(\frac{(2i-1)\pi}{2d}(k-\varepsilon_0)\right)\right)$$

$$=\frac{d}{2^{2t+1}}\left(C_{2t+1}^r+\sum_{i>0}(-1)^i(C_{2t+1}^{t-\varepsilon_0-i\lambda}+C_{2t+1}^{t+\varepsilon_0+1-i\lambda})\right) \qquad \text{by (2-b)}$$

$$=\frac{d}{2^{2t+1}}\left(C_m^r+\sum_{i>0}(-1)^i(C_m^{r-i\lambda}+C_m^{s-i\lambda})\right)$$

$$=\frac{\lambda}{2^{2t+3}}\cdot|\mathfrak{D}(m,\varepsilon,\delta)|.$$

i.e.

$$|\mathfrak{D}(m,\varepsilon,\delta)|=\frac{2^{m+2}}{\lambda}\cdot\sum_{1\le i\le d}\cos(\varepsilon(2i-1)\pi/(2\lambda))\cdot\cos^m((2i-1)\pi/(2\lambda)).$$

If $\delta=2d$, then

$$\sum_{1\le i\le d}\cos(\varepsilon((2i-1)\pi/(2\delta+2))\cdot\cos^{2t+1}((2i-1)\pi/(2\delta+2))$$

$$=\frac{1}{2^{2t+1}}\cdot\sum_{0\le k\le t}C_{2t+1}^{t-k}\cdot\sum_{1\le i\le d}\left(\cos\left(\frac{(2i-1)\pi}{(2\delta+2)}((2k+1)+\varepsilon)\right)+\cos\left(\frac{(2i-1)\pi}{(2\delta+2)}((2k+1)-\varepsilon)\right)\right)$$

$$=\frac{1}{2^{2t+1}}\cdot\sum_{0\le k\le t}C_{2t+1}^{t-k}\cdot\sum_{1\le i\le d}\left(\cos\left(\frac{(2i-1)\pi}{(2d+1)}((k+1)+\varepsilon_0)\right)+\cos\left(\frac{(2i-1)\pi}{(2d+1)}(k-\varepsilon_0)\right)\right)$$

$$=\frac{1}{2^{2t+1}}\left(d\cdot C_{2t+1}^{t-\varepsilon_0}+d\cdot\sum_{i>0}(-1)^i(C_{2t+1}^{t-\varepsilon_0-i\lambda}+C_{2t+1}^{t+\varepsilon_0+1-i\lambda})\right)$$

$$+\frac{1}{2^{2t+1}}\left(\frac{1}{2}\sum_{((k+1)+\varepsilon_0)\nmid\lambda}(-1)^{k+\varepsilon_0}C_{2t+1}^{t-k}+\frac{1}{2}\sum_{(k-\varepsilon_0)\nmid\lambda}(-1)^{k-\varepsilon_0-1}C_{2t+1}^{t-k}\right) \quad \text{by (2-c)}$$

$$=\frac{1}{2^{2t+1}}\left(d\cdot C_{2t+1}^{t-\varepsilon_0}+d\cdot\sum_{i>0}(-1)^i(C_{2t+1}^{t-\varepsilon_0-i\lambda}+C_{2t+1}^{t+\varepsilon_0+1-i\lambda})\right)$$

$$+\frac{1}{2^{2t+1}}\left(\frac{1}{2}\sum_{\substack{((k+1)+\varepsilon_0)\nmid\lambda \\ (k-\varepsilon_0)\mid\lambda}}(-1)^{k+\varepsilon_0}C_{2t+1}^{t-k}+\frac{1}{2}\sum_{\substack{(k-\varepsilon_0)\nmid\lambda \\ ((k+1)+\varepsilon_0)\mid\lambda}}(-1)^{k-\varepsilon_0-1}C_{2t+1}^{t-k}\right)$$

$$=\frac{1}{2^{2t+1}}\left(\left(d+\frac{1}{2}\right)\cdot C_{2t+1}^{t-\varepsilon_0}+\left(d+\frac{1}{2}\right)\cdot\sum_{i>0}(-1)^i(C_{2t+1}^{t-\varepsilon_0-i\lambda}+C_{2t+1}^{t+\varepsilon_0+1-i\lambda})\right)$$

$$=\frac{\lambda}{2^{2t+2}}\left(C_m^r+\sum_{i>0}(-1)^i(C_m^{r-i\lambda}+C_m^{s-i\lambda})\right)$$

$$=\frac{\lambda}{2^{2t+3}}\cdot|\mathfrak{D}(m,\varepsilon,\delta)|.$$

i.e.

$$|\mathfrak{D}(m,\varepsilon,\delta)|=\frac{2^{m+2}}{\lambda}\cdot\sum_{1\leq i\leq d}\cos(\varepsilon(2i-1)\pi/(2\lambda))\cdot\cos^m((2i-1)\pi/(2\lambda)).$$

The proof for the case $m$ odd has been completed, the proof for the case $m$ even is similar provided applying formula (2-a), so which is omitted. $\square$

*Proof of Theorem 10:* At first, it is noticed that $\varepsilon$ and $m$ should be same parity. So, by (2.42), If $m$ is odd, $\delta=2d-1$, then

$$|\mathfrak{S}(m,\delta)|=\sum_{\varepsilon\leq\delta}|\mathfrak{D}(m,\varepsilon,\delta)|=\sum_{1\leq i\leq d}|\mathfrak{D}(m,2i-1,\delta)|$$

$$=\frac{2^{m+2}}{\lambda}\cdot\sum_{1\leq i\leq d}\sum_{j=1}^d\cos((2j-1)(2i-1)\pi/(2\lambda))\cdot\cos^m((2i-1)\pi/(2\lambda))$$

$$=\frac{2^{m+1}}{\lambda}\cdot\sum_{1\leq i\leq d}\frac{(-1)^{i-1}}{\sin((2i-1)\pi/(2\lambda))}\cdot\cos^m((2i-1)\pi/(2\lambda)) \quad \text{by (2-f)}$$

If $m$ is odd, $\delta=2d$, then

$$|\mathfrak{S}(m,\delta)|=\sum_{\varepsilon\leq\delta}|\mathfrak{D}(m,\varepsilon,\delta)|=\sum_{1\leq i\leq d}|\mathfrak{D}(m,2i-1,\delta)|$$

$$=\frac{2^{m+2}}{\lambda}\cdot\sum_{1\leq i\leq d}\sum_{j=1}^d\cos((2j-1)(2i-1)\pi/(2\lambda))\cdot\cos^m((2i-1)\pi/(2\lambda))$$

$$=\frac{2^{m+1}}{\lambda}\cdot\sum_{1\leq i\leq d}(-1)^{i-1}\cot\left((2i-1)\pi/(2\lambda)\right)\cdot\cos^m((2i-1)\pi/(2\lambda)) \quad \text{by (2-g)}$$

If $m$ is even, $\delta=2d-1$, then

$$|\mathfrak{S}(m,\delta)|=\sum_{\varepsilon\leq\delta}|\mathfrak{D}(m,\varepsilon,\delta)|$$

$$= \frac{2^{m+1}}{\lambda} \cdot \sum_{1 \leq i \leq d} \sum_{j=0}^{d-1} (1+\text{sgn}(j)) \cdot \cos((2j)(2i-1)\pi/(2\lambda)) \cdot \cos^m((2i-1)\pi/(2\lambda))$$

$$= \frac{2^{m+1}}{\lambda} \cdot \sum_{1 \leq i \leq d} (-1)^{i-1} \cot\left((2i-1)\pi/(2\lambda)\right) \cdot \cos^m((2i-1)\pi/(2\lambda)) \qquad \text{by (2-h)}$$

If $m$ is even, $\delta = 2d$, then

$$|\mathfrak{S}(m,\delta)| = \sum_{\varepsilon \leq \delta} |\mathfrak{D}(m,\varepsilon,\delta)|$$

$$= \frac{2^{m+1}}{\lambda} \cdot \sum_{1 \leq i \leq d} \sum_{j=0}^{d} (1+\text{sgn}(j)) \cos((2j)(2i-1)\pi/(2\lambda)) \cdot \cos^m((2i-1)\pi/(2\lambda))$$

$$= \frac{2^{m+1}}{\lambda} \cdot \sum_{1 \leq i \leq d} \frac{(-1)^{i-1}}{\sin((2i-1)\pi/(2\lambda))} \cdot \cos^m((2i-1)\pi/(2\lambda)) \ . \qquad \text{by (2-i)}$$

□

In the following, we give another proof for Theorem 9, which is similar to the one of Theorem 3, by applying the induction and a generalized Newton formula.

It is clear that the original Newton's formula now is no longer applicable here. So, the first is to develop the formula.

Suppose that $S = \{1, 2, \cdots, d\}$, denoted by $S^{(k)}$ the set of all $k$-subsets of $S$. Suppose that $x_1, x_2, \cdots, x_d$ are $d$ variables, for $A \in S^{(k)}$, denoted by $x_A = \sum_{i \in A} x_i$, and $\hat{x}_A = \prod_{i \in A} x_i$ respectively. Suppose that $f(x)$ is a monic polynomial of degree $d$, $f(x) = \sum_{0 \leq i \leq d} a_i x^{d-i}$, ($a_0 = 1$), and $x_1, x_2, \cdots, x_d$ are the $d$ zeros of $f(x)$, $c = (c_1, c_2, \cdots, c_d)$, $c_i$, $i = 1, 2, \cdots, d$ are $d$ constants. For a non-negative integer $k$, denoted by $\rho_k(c) = \sum_{1 \leq i \leq d} c_i x_i^k$, and $\tau_k(c) = (-1)^k \cdot \sum_{A \in S^{(k)}} c_A \cdot \hat{x}_A$, then there is

**Theorem A**. (generalization of Newton's formula)

$$\begin{cases} \sum_{0 \leq k \leq r-1} \rho_{r-k}(c) \cdot a_k + \tau_r(c) = 0, & \text{for } r \leq d, \qquad (2.44) \\ \sum_{0 \leq k \leq d} \rho_{r-k}(c) \cdot a_k = 0, & \text{for } r > d. \qquad (2.45) \end{cases}$$

*Proof.* If $r \leq d$, then

$$\sum_{0 \leq k \leq r-1} \rho_{r-k}(c) \cdot a_k = \rho_r(c) + \sum_{1 \leq k \leq r-1} \left( \sum_{1 \leq i \leq d} c_i x_i^{r-k} \right) \cdot \left( \sum_{A \in S^{(k)}} (-1)^k \hat{x}_A \right)$$

$$= \rho_r(c) + \sum_{1 \le k \le r-1} \left( (-1)^k \sum_{A \in S^{(k)}} \left( \sum_{i \in A} c_i x_i^{r-k} \right) \cdot \hat{x}_A + (-1)^k \sum_{B \in S^{(k+1)}} \left( \sum_{i \in B} c_i x_i^{r-k-1} \right) \cdot \hat{x}_B \right)$$

$$= (-1)^{r-1} \sum_{B \in S^{(r)}} c_B \cdot \hat{x}_B$$

$$= -\tau_r(c).$$

If $r > d$, then

$$\sum_{0 \le k \le d} \rho_{r-k}(c) \cdot a_k = \sum_{0 \le k \le d} a_k \cdot \sum_{1 \le i \le d} c_i x_i^{r-k}$$

$$= \sum_{1 \le i \le d} c_i \cdot \left( \sum_{0 \le k \le d} a_k x_i^{r-k} \right)$$

$$= 0.$$

□

Clearly, in the case $c = (1, 1, \cdots, 1)$ the formula is just right the original Newton's formula.

Simply write $z_i = (2i-1)\pi / (2\delta + 2)$, for our application here, take $f(x) = x^d \cdot p_\delta(x^{-1})$, $x_i = \vartheta_i^{-1}$, $c_i = \cos(\varepsilon \cdot z_i)$, $1 \le i \le d$, and in this application, we specially write $a_k$, $\rho_k(c)$ and $\tau_k(c)$ as $a_{\delta,k}$, $\rho_k$ and $\tau_{\delta,k}(\varepsilon)$ respectively to keep consistency in the context.

By the formula (2-c), it is easy to follow that

$$\sum_{i \in S} \left( 4\cos^2(z_i) \cdot \cos((\varepsilon - 2)z_i) \right) = \begin{cases} \lambda, & \text{if } \varepsilon = 2, \\ \frac{1}{2}\lambda, & \text{if } \varepsilon = 0, \text{ or } 4, \\ 0, & \text{if } \varepsilon > 4, \text{ even.} \end{cases} \quad (2.46)$$

Let $\mu = \frac{1}{2}$, or $1$, or $0$, as $\varepsilon = 0, 4$, or $2$, or else. With (2.46), it has

**Lemma 11.** Suppose that $\varepsilon$ is an even number, then there is

$$\tau_{\delta,r}(\varepsilon) - \tau_{\delta,r+1}(\varepsilon - 2) + 2 \cdot \tau_{\delta,r}(\varepsilon - 2) + \tau_{\delta,r}(\varepsilon - 4) = \mu \cdot (\delta + 1) \cdot a_{\delta,r}. \quad (2.47)$$

*Proof.*

$$\tau_{\delta,r}(\varepsilon) + \tau_{\delta,r}(\varepsilon - 4) = (-1)^r \cdot \sum_{A \in S^{(r)}} \sum_{i \in A} \left( \cos(\varepsilon z_i) + \cos((\varepsilon - 4)z_i) \right) \cdot \hat{x}_A$$

$$= (-1)^r \sum_{A \in S^{(r)}} \sum_{i \in A} \left( 2\cos(2z_i) \cdot \cos((\varepsilon - 2)z_i) \right) \cdot \hat{x}_A$$

$$= (-1)^r \sum_{A \in S^{(r)}} \sum_{i \in A} \left( 4\cos^2(z_i) \cdot \cos((\varepsilon - 2)z_i) \right) \cdot \hat{x}_A - (-1)^r \cdot 2 \sum_{A \in S^{(r)}} \sum_{i \in A} \cos((\varepsilon - 2)z_i) \cdot \hat{x}_A$$

$$= (-1)^r \sum_{A \in S^{(r)}} \sum_{i \in A} \left( 4\cos^2(z_i) \cdot \cos((\varepsilon - 2)z_i) \right) \cdot \hat{x}_A - 2 \cdot \tau_{\delta,r}(\varepsilon - 2)$$

$$= (-1)^r \sum_{i \in S} \left(4\cos^2(z_i) \cdot \cos((\varepsilon-2)z_i)\right) \cdot \sum_{A \in S^{(r)}} \hat{x}_A - 2 \cdot \tau_{\delta,r}(\varepsilon-2)$$

$$+ (-1)^{r+1} \sum_{B \in S^{(r+1)}} \left(\sum_{i \in B} \cos((\varepsilon-2)z_i)\right) \cdot \hat{x}_B$$

$$= \tau_{\delta,r+1}(\varepsilon-2) - 2 \cdot \tau_{\delta,r}(\varepsilon-2) + \mu \cdot (\delta+1) \cdot a_{\delta,r}. \qquad \text{by (2.46)}$$

□

It is clear that $\tau_{\delta,r}(0) = r \cdot a_{\delta,r}$, so, with (2.25), (2.47), and correlation $a_{\delta,r} = a_{\delta-1,r} - a_{\delta-2,r-1}$, take the induction on $\varepsilon$, it may be followed directly that

**Corollary 4.**

$$\frac{\tau_{\delta,r}(\varepsilon)}{\delta+1} = \frac{\tau_{\delta-1,r}(\varepsilon)}{\delta} - \frac{\tau_{\delta-2,r-1}(\varepsilon)}{\delta-1}. \qquad (\varepsilon \text{ even}) \qquad (2.48)$$

And,

**Corollary 5.** For two positive integers $r$ and $\delta$, denoted by $\tilde{\omega}_r(\varepsilon,\delta) = |\mathfrak{D}(2r,\varepsilon,\delta)|$, and let $\tilde{\omega}_0(\varepsilon,\delta) = \frac{2^{1+\text{sgn}(\varepsilon)}}{\delta+1} \tau_{\delta,r}(\varepsilon)$ if $r \leq d$, or $0$ otherwise, then it has

$$\sum_{0 \leq i < r} \tilde{\omega}_{r-i}(\varepsilon,\delta) \cdot a_{\delta,i} + \tilde{\omega}_0(\varepsilon,\delta) = 0. \qquad (2.49)$$

Corollary 5 may be proved with (2.48) in a similar way as one of Corollary 3, but it should be noticed that there are a little difference that here is $\tilde{\omega}_r(\delta) = \tilde{\omega}_r(\delta-1) = \cdots = \tilde{\omega}_r(r')$, $r' = (m+\varepsilon)/2$, so there are two special cases: $m = 2d, \delta = 2d-1, \varepsilon = 2d-2$, and $m = 2d, \delta = 2d, \varepsilon = 2d$, which need to be verified separately.

It is easy to obtain that

$$\sum_{1 \leq i \leq d} \cos(\varepsilon z_i) = \begin{cases} 0, & \text{if } \delta = 2d-1, \varepsilon = 2k(>0). \\ \dfrac{(-1)^{k-1}}{2}, & \text{if } \delta = 2d, \varepsilon = 2k(>0). \end{cases}$$

And by (2.19), it can be known that $a_{2d,d} = (-1)^d \cdot (2d+1)$, hence

$$\frac{4}{\delta+1} \cdot \tau_{\delta,d}(\varepsilon) = \frac{4}{\delta+1} \cdot \sum_{1 \leq i \leq d} \cos(\varepsilon z_i) \cdot a_{\delta,k} = \begin{cases} 0, & \text{if } \delta = 2d-1, \varepsilon = 2k(>0). \\ -2, & \text{if } \delta = 2d, \varepsilon = 2d(>0). \end{cases}$$

Thus, (2.49) is also held for the two special cases.

*The Second Proof for Theorem 9:*

As the proof of Theorem 3, we take the induction on $m$.

At first, we assume that $m(=2r)$ is even. By (2.46) we know that (2.42) stand for $m=2$.

Generally, for $r \leq d$, by induction, with (2.44) and (2.49), it has

$$\rho_r = -\left(\sum_{1\leq k \leq r-1} \rho_{r-k} a_{\delta,k} + \tau_{\delta,r}\right) = -\left(\sum_{1\leq k \leq r-1} \frac{\delta+1}{2^{1+\mathrm{sgn}(\varepsilon)}} \cdot \tilde{\omega}_{r-k} \cdot a_{\delta,k} + \tau_{\delta,r}\right)$$

$$= -\frac{\delta+1}{2^{1+\mathrm{sgn}(\varepsilon)}} \cdot \left(\sum_{1\leq k \leq r-1} \tilde{\omega}_{r-k} a_{\delta,k} + \frac{2^{1+\mathrm{sgn}(\varepsilon)}}{\delta+1} \tau_{\delta,r}\right)$$

$$= \frac{\delta+1}{2^{1+\mathrm{sgn}(\varepsilon)}} \tilde{\omega}_r$$

i.e.

$$\tilde{\omega}_r = \frac{2^{1+\mathrm{sgn}(\varepsilon)}}{(\delta+1)} \rho_r.$$

For $r > d$, it may be shown by (2.45) and the induction.

For the case $m$ is odd, and so $\varepsilon$ is odd, the difference from the one of the case $m$ even is that here it should be taken $c_i = \cos(\varepsilon \cdot z_i)/\cos(z_i)$, however, we know $\cos((2k+1)z)/\cos(z)$ is a linear sum of $\cos(2iz)$, $i = 0,1,\cdots,k$. Thus, the approach above is also available for the case $m$ odd, the details are omitted. □

Finally, the author would like to take the opportunity to acknowledge arXiv and the other websites publishing the paper, and thank all the readers.

# Appendix    A Reference Code

## to find the values of $\zeta(m,\delta)$ by counting and by formula respectively

```c
#include <stdio.h>
#include <math.h>

main()
{
int i,k,id,num,delta;
for(num=18;num<32;num+=2)
{
for(delta=1;delta<=12;delta++)
{
unsigned int   n=(1<<num);
unsigned int x,sum=0;
for(k=1;k<n;k+=2)
{
int t=0;
int a[32]={0};
x=k;
id=(x&1);
for(i=0;i<num;i++)
{
if((x&1)==id){
a[t]++;
}
else{
t=t+1;
a[t]++;
id=(id^1);
}
x>>=1;
}
int y=0;
x=0;
for(i=0;i<num;i+=2)
{
x+=a[i];
y+=a[i+1];
}
int h;
if(x>=y)h=x-y;
else h=y-x;
```

```c
if(h<=1){
int v,s=0;
t=(-1);
for(i=0;i<num;i++)
{
t=t*(-1);
s+=(a[i]*t);
if(s>=0)v=s;
else v=(-s);
if(v>delta)break;
}
if(v<=delta){
sum++;
}
}
}
}
double z=0,u,c;
int d=floor((delta+1)/2);
k=(num+1)>>1;
for(i=0;i<d;i++){
c=cos(3.14159265359*(0.5+i)/(delta+1));          //pai=3.14159265359
u=4*c*c;
z+=pow(u,k);
}
z=z/(delta+1);
printf("%d     %d     %d     %f\n",num, delta, sum,z);
}
}
}
```

/*********************************************************************/